\documentclass[10pt]{article}
\usepackage{amssymb, amsfonts, amsthm, mathrsfs}
\usepackage[super]{nth}
\usepackage[pagebackref,colorlinks=true]{hyperref}
\usepackage{graphicx, xcolor}
\usepackage[all]{xy}
\usepackage[leqno]{amsmath}
\usepackage{titlesec}

\theoremstyle{definition}
\newtheorem{thm}{Theorem}[section]

\newtheorem{cor}[thm]{Corollary}
\newtheorem{prop}[thm]{Proposition}

\titleformat{\section}{\Large\bfseries}{\thesection}{.5em}{}

\makeatletter

\def\id{\mathrm{id}}

\def\tr{\mathrm{Tr}}

\def\ccc{\mathbb{C}}
\def\C{\mathbb{C}}

\def\rr{\mathbb{R}}
\def\R{\mathbb{R}}
\def\pp{\mathbb{P}}

\def\clo{\mathcal{O}}

\def\frI{\mathfrak{I}}

\def\pt{\partial}
\def\p{\partial}
\def\bpt{\bar{\pt}}

\def\ddb{\pt\bpt}

\def\ud{\mathrm{d}}

\def\bw{\bar{w}}
\def\bz{\bar{z}}

\def\bzeta{\bar{\zeta}}

\def\vol{\mathrm{vol}}

\def\U{\mathrm{U}}
\def\SU{\mathrm{SU}}

\def\spann{\mathrm{Span}}

\def\ind{\mathrm{Ind}}
\makeatother

\begin{document}

\title{\textbf{Generalized Calabi-Gray Geometry and Heterotic Superstrings}}
\author{Teng Fei}

\date{}
\maketitle{}

\section{Introduction}

The Hull-Strominger system describes the geometry of compactification of heterotic superstrings with torsion to 4-dimensional Minkowski spacetime. In physics, one first writes down the bosonic part of the heterotic superstring action and then computes its critical points to get the equations of motion. Preferably, physicists would like to solve the equations of motion with unbroken supersymmetry, which is encoded in the Killing spinor equation. There were arguments that supersymmetry implies the equations of motion, however this feature has not been clarified mathematically.
Additionally, for the theory to be consistent, the Bianchi identity must be satisfied.

The fields and geometric objects involved in the aforementioned action are a 10-dimensional manifold $M$ with a Lorentzian metric $g$, a vector bundle $E$ over $M$ with curvature tensor $F$, a closed 3-form $H$ known as the NS field strength, and a scalar dilaton function $\phi$. For practical reason, we expect the theory to be compactified to a 4-dimensional spacetime. That is, $M$ is a product of a 4-dimensional spacetime $S$ and a compact 6-dimensional manifold $X$ known as the internal space. Moreover, all the fields have only $X$-dependence.

Such a reduction was first successfully carried out by Candelas-Horowitz-Strominger-Witten in \cite{candelas1985}, where they take $(M,g)$ to be a metric product of the Minkowski spacetime and a Riemannian 6-manifold $X$. By embedding the spin connection into the gauge connection, one can cancel the anomaly as required by Green-Schwarz \cite{green1984}. Moreover the supersymmetry condition is equivalent to that $X$ is a K\"ahler manifold with Ricci-flat metric, i.e., a Calabi-Yau space.

To allow for compactifications with torsion, one needs to replace the metric product by a warped product. This idea was implemented by Hull \cite{hull1986c, hull1986b} and Strominger \cite{strominger1986} independently, where they showed that, among many things else, supersymmetry implies that $X$ has to be a complex manifold with holomorphically trivial canonical bundle. We shall refer the system of PDE's derived in Hull and Strominger's papers as the \emph{Hull-Strominger system}. It should be emphasized that the Hull-Strominger system incorporates only the Bianchi identity and $\mathcal{N}=1$ supersymmetry, which has no implication on equations of motion. Besides the original derivations, one may consult \cite{wu2006} for a more mathematical treatment.

In this note, we regard the Hull-Strominger system as a pure mathematical problem, which can be formulated as follows after Li-Yau \cite{li2005}:
\begin{eqnarray}
\label{hym}F\wedge\omega^2=0,\quad F^{0,2}=F^{2,0}=0,\\
\label{ac}i\pt\bpt\omega=\frac{\alpha'}{4}(\tr(R\wedge R)-\tr(F\wedge F)),\\
\label{cb}\ud(\|\Omega\|_\omega\cdot\omega^2)=0.
\end{eqnarray}
Here $X$ is a complex manifold equipped with a Hermitian metric $\omega$. In addition, $X$ has holomorphically trivial canonical bundle, which is trivialized by a nowhere vanishing holomorphic volume form $\Omega$. Besides that, $E\to X$ is a holomorphic vector bundle equipped with a Hermitian metric $H$. Denoted by $R$ and $F$ the endomorphism valued curvature tensors of the holomorphic tangent bundle $T^{1,0}X$ and $E$ respectively. Finally $\alpha'$ is a coupling constant, which is always positive in this paper. Fixing $X$ and $\Omega$, the goal is to find metrics $\omega$ and $H$ such that (\ref{hym})-(\ref{cb}) hold.

Suppose $\omega$ is a K\"ahler metric, i.e. $\ud\omega=0$, then (\ref{cb}) implies that $\|\Omega\|_\omega$ is a constant and the metric is Ricci-flat. Due to Yau's solution \cite{yau1978} to the Calabi conjecture, such $\omega$ exists uniquely in each K\"ahler class. Moreover, by taking $E=T^{1,0}X$ and setting $F=R$, (\ref{hym}) and (\ref{ac}) are solved simultaneously and we get exactly the torsion-free compactification proposed in \cite{candelas1985}.

The novel feature for the Hull-Strominger system is that it allows potentially for non-K\"ahler backgrounds: the K\"ahler condition $\ud\omega=0$ is relaxed to (\ref{cb}). If we consider the conformally changed metric $\tilde{\omega}=\|\Omega\|^{1/2}_\omega\cdot\omega$, then (\ref{cb}) can be phrased as $\ud(\tilde{\omega}^2)=0$, which is exactly the definition of $\tilde{\omega}$ being a balanced metric. Balanced metrics have been extensively studied in mathematics and their existence are well-understood \cite{michelsohn1982}. For this reason, we call (\ref{cb}) the \emph{conformally balanced equation}, its solvability is equivalent to the existence of balanced metrics on $X$.

Equation (\ref{hym}) is the famous \emph{Hermitian Yang-Mills equation}. As we can replace $\omega$ by $\tilde{\omega}$ in (\ref{hym}), Li-Yau's Gauduchon version \cite{li1987} of Donaldson-Uhlenbeck-Yau theorem applies, therefore (\ref{hym}) is solvable if and only if $E$ is polystable with respect to the polarization $\tilde{\omega}$.

Geometry of the manifold $X$ and of the bundle $E$ are coupled together via the Bianchi identity (\ref{ac}), to which we also refer as the \emph{anomaly cancellation equation}. This is an equation of 4-forms on $X$, which is in general very hard to tackle. Taking cohomology classes of both sides, we know that $c_2(E)$, the second Chern class of $E$, is the same as the $c_2(X)$, provided $E$ has trivial $c_1$.

There are many aspects of the Hull-Strominger system, such as construction of solutions, geometry of moduli spaces, viewpoints from generalized geometry and so on. As there is already an excellent introductory text \cite{garciafernandez2016b} to many of these themes, especially to the second and the third aspects listed above, we will focus only on construction of solutions in this note.

This paper is organized as follows. In Section 2, we give a brief review of known solutions to the Hull-Strominger system in literature. In Section 3 we present the geometric construction of generalized Calabi-Gray manifolds and solutions to the Hull-Strominger system on them. Section 4 will be devoted to Anomaly flow and its reduction on generalized Calabi-Gray manifolds. Due to its expository nature, most parts of this paper have appeared in the literature before and there is no intention for originality.

\section{A Tour of Known Solutions}

Because the Hull-Strominger system comes from compactifications of heterotic superstrings, we are mostly interested in solving the Hull-Strominger system on compact manifolds. Nevertheless it is also useful to study solutions on noncompact manifolds (also known as local models) such as \cite{fu2009, carlevaro2010, halmagyi2016, fei2017b}, as they may provide important information in gluing constructions and in understanding the string landscape \cite{reid1987, fu2012, chuan2011}.

As we have seen in the introduction, K\"ahler Ricci-flat metrics are solutions to the Hull-Strominger system. The simplest way to obtain non-K\"ahler solutions is to deform the K\"ahler ones. This idea was developed in \cite{strominger1986, li2005, andreas2012, andreas2012b} and indeed one may obtain non-K\"ahler solutions to the Hull-Strominger system on any K\"ahler Calabi-Yau manifold with suitable holomorphic vector bundles.

However, the novel feature for the Hull-Strominger system is that it allows for non-K\"ahler backgrounds. Therefore it is of great interest to find essentially non-K\"ahler solutions, that is, to construct solutions to the Hull-Strominger system on complex manifolds with no K\"ahler structure at all. The first such example was found by Fu-Yau \cite{fu2008}. More precisely, Fu and Yau studied the Hull-Strominger system on a class of toric fibrations over K3 surfaces constructed in \cite{goldstein2004} dating back to Calabi-Eckmann \cite{calabi1953b}, and showed that, by choosing a suitable ansatze, the Hull-Strominger system reduces to a complex Monge-Amp\`ere type equation of a scalar function on K3 surfaces. By hard analysis, Fu-Yau proved the existence of solution to this fully nonlinear PDE, hence constructed solutions to the Hull-Strominger system on a class of non-K\"ahler Calabi-Yau 3-folds. The Fu-Yau equation has rich content in both mathematics and physics, therefore it has inspired many later developments including \cite{becker2006, fu2007, becker2009, phong2016a, phong2017c, phong2018, chu2018}.

Another class of non-K\"ahler solutions to the Hull-Strominger system are constructed on Lie groups and their quotient by discrete subgroups, including \cite{fernandez2009, grantcharov2011, fernandez2014, ugarte2014, ugarte2015, fei2015, otal2017}. The main idea of this construction is to make everything left-invariant, therefore the Hull-Strominger system reduces to a system of algebraic equations on corresponding Lie algebras.

For all the compact solutions mentioned above, including the K\"ahler ones, there are only finitely many sets of Hodge numbers for the internal space $X$ (there may be infinitely many different topological types in the Goldstein-Prokushkin/Fu-Yau construction). In fact in K\"ahler geometry, it is a well-known conjecture that there are only finitely many deformation types, thus finitely many sets of Hodge numbers, of compact K\"ahler Calabi-Yau manifolds in each dimension. Since the Hull-Strominger system can be viewed as a generalization of K\"ahler Calabi-Yau's, it is natural to ask the same question for compact complex manifolds admitting solutions to the Hull-Strominger system. One of the main results in this article is to give a negative answer to this question. More accurately, we present an explicit construction of solutions to the Hull-Strominger system on generalized Calabi-Gray manifolds, which are compact non-K\"ahler Calabi-Yau 3-folds with infinitely many topological types and sets of Hodge numbers.

\section{The Geometry of Generalized Calabi-Gray Manifolds}

\subsection*{Construction}

In \cite{calabi1958}, Calabi constructed a class of complex non-K\"ahler manifolds with natural $\SU(3)$ structures to show that the first Chern class $c_1$ is not a smooth invariant. Calabi's construction was later generalized by Gray \cite{gray1969} in the setting of manifolds with vector cross products. This famous Calabi-Gray construction relies on the exotic 3-fold vector cross product on 7-manifolds with $\mathrm{G}_2$-structures. It was later realized \cite{fei2015b, fei2016b} that Calabi-Gray manifolds have holomorphically trivial canonical bundles and they admit balanced metrics. As a consequence, they provide fantastic candidates for solutions to the Hull-Strominger system.

As observed in \cite{fei2016}, the classical Calabi-Gray construction is related to twistor spaces of hyperk\"ahler 4-manifolds, which further allows us to give a small yet essential generalization of Calabi and Gray's work. Following \cite{fei2017}, we shall describe the construction of generalized Calabi-Gray manifolds as follows.

Let $(M,g)$ be a compact hyperk\"ahler 4-manifold, that is, $(M,g)$ is a Riemannian 4-manifold equipped with three compatible complex structures $I$, $J$ and $K$ satisfying $IJK=-\id$. It is well-known that $(M,g)$ is either a flat 4-torus or a K3 surface with a Calabi-Yau metric. Denote by $\omega_I$, $\omega_J$ and $\omega_K$ the corresponding K\"ahler forms, defined as usual by $\omega_I(v_1,v_2)=g(Iv_1,v_2)$, $\omega_J(v_1,v_2)=g(Jv_1,v_2)$ and $\omega_K(v_1,v_2) = g(Kv_1,v_2)$. In fact, for any real numbers $\alpha$, $\beta$ and $\gamma$ such that
\[\alpha^2+\beta^2+\gamma^2=1,\]
the combination $\alpha I+\beta J+\gamma K$ is also a compatible complex structure, whose associated K\"ahler form is given by $\alpha\omega_I+\beta\omega_J+\gamma\omega_K$.

Following \cite{hitchin1987}, let us first describe the twistor space $Z$ of $M$. By stereographic projection, we may parameterize $S^2=\{(\alpha,\beta,\gamma)\in\rr^3:\alpha^2+\beta^2+\gamma^2=1\}$ by $\zeta\in\ccc\pp^1$ via
\begin{equation}\label{stereographic}
(\alpha,\beta,\gamma)= \left(\frac{1-|\zeta|^2}{1+|\zeta|^2},\frac{\zeta+\bzeta}{1+|\zeta|^2},\frac{i(\bzeta-\zeta)}{1+|\zeta|^2}\right).
\end{equation}
We always fix a round (Fubini-Study) metric
\begin{equation}\label{fs}
\omega_{\textrm{FS}}=\frac{2i\ud\zeta\wedge\ud\bzeta}{(1+|\zeta|^2)^2}
\end{equation}
on $S^2\cong\ccc\pp^1$.

The twistor space $Z$ of $M$ is defined to be the manifold $Z=\ccc\pp^1\times M$ with the tautological almost complex structure $\frI$ given by
\begin{equation}\label{hypertwistor}
\frI_{(\zeta,x)}=j_\zeta\oplus(\alpha(\zeta) I_x+\beta(\zeta) J_x+\gamma(\zeta) K_x)
\end{equation}
for $(\zeta,x)\in\ccc\pp^1\times M$, where $j$ is the standard complex structure on $\ccc\pp^1$ with holomorphic coordinate $\zeta$.

In fact, we have
\begin{thm}\cite{atiyah1978, hitchin1981, hitchin1987, muskarov1989}
\begin{enumerate}
\item $\frI$ is an integrable complex structure, making $(Z,\frI)$ a compact non-K\"ahler 3-fold whose natural product metric is balanced.
\item The natural projection $\pi:Z=\ccc\pp^1\times M\to\ccc\pp^1$ is holomorphic.
\item Let $\wedge^2\Omega_{Z/\ccc\pp^1}$ be the determinant line bundle of the relative cotangent bundle associated to the holomorphic fibration $\pi$. Then the line bundle $\wedge^2\Omega_{Z/\ccc\pp^1}\otimes\pi^*\clo(2)$ on $Z$ has a global section which defines a holomorphic symplectic form on each fiber of $\pi$.
\end{enumerate}
\end{thm}

Now let $\Sigma$ be a compact Riemann surface of genus $g$ and let $\varphi:\Sigma\to\ccc\pp^1$ be a nonconstant holomorphic map. We may treat $\varphi$ either as a nonconstant meromorphic function $\zeta$ on $\Sigma$ or, by identifying $\ccc\pp^1$ with the unit sphere in $\rr^3$, as a map $\varphi=(\alpha,\beta,\gamma)$ into $\rr^3$. These two viewpoints are related by the stereographic projection formula (\ref{stereographic}). By pulling back the holomorphic fibration $\pi:Z\to\ccc\pp^1$ over $\varphi:\Sigma\to\ccc\pp^1$, we get a holomorphic fibration $p:X=\varphi^*Z\to\Sigma$. As a complex manifold, $X$ is topologically $\Sigma\times M$ with a twisted complex structure $J_0 = j_\Sigma \oplus(\alpha I_x+\beta J_x+\gamma K_x)$. We have
\begin{prop}\cite{fei2015b, fei2016, fei2016b}
\begin{enumerate}
\item $X$ has holomorphically trivial canonical bundle if and only if
\begin{equation}\label{spinor}
\varphi^*\clo(2)\cong K_\Sigma,
\end{equation}
where $K_\Sigma$ is the canonical bundle of $\Sigma$.
\item Under (a), $X$ is non-K\"ahler with balanced metrics.
\end{enumerate}
\end{prop}

From now on, we shall call condition (\ref{spinor}) plus that $\varphi$ is not the constant map the ``pullback condition''. Assuming the pullback condition, it is clear that $S=\varphi^*\clo(1)$ is a square root of $K_\Sigma$, which is known as a theta characteristic in algebraic geometry, or alternatively a spin structure according to Atiyah \cite{atiyah1971}. Furthermore, the linear system associated to the line bundle $S$ is basepoint-free.

Conversely, if we start with a basepoint-free theta characteristic $S$ on $\Sigma$, we may choose $s_1,s_2\in H^0(\Sigma,S)$ such that $s_1$ and $s_2$ have no common zeroes, then $\zeta=s_1/s_2$ is a meromorphic function on $\Sigma$. Moreover, $\zeta$ defines a holomorphic map $\varphi:\Sigma\to\ccc\pp^1$ such that the pullback condition holds.

Therefore in this expose, we shall call a pair $(\Sigma,\varphi)$ such that the pullback condition holds a \emph{vanishing spinorial pair}. The generalized Calabi-Gray construction says that from any vanishing spinorial pair, one can construct a compact non-K\"ahler Calabi-Yau 3-fold $X$ with balanced metrics. It is trivial to remark that if $(\Sigma,\varphi)$ is a vanishing spinorial pair and $\tau:\widetilde{\Sigma}\to\Sigma$ an unramified covering map, then $(\widetilde{\Sigma},\varphi\circ\tau)$ is also a vanishing spinorial pair. We also say such an $X$ is a \emph{genus $g$ generalized Calabi-Gray manifold}, where $g$ is the genus of $\Sigma$.

\subsection*{Hermitian Metrics}

As a smooth manifold $X$ is a product of $\Sigma$ and $M$, and $M$ is equipped with a fixed hyperk\"ahler metric, therefore it is very useful to construct a canonical metric $\hat{\omega}$ on $\Sigma$ from the vanishing spinorial pair $(\Sigma,\varphi)$. A guiding example is the classical Calabi-Gray construction, where $\Sigma$ is an immersed minimal surface in a flat 3-torus and $\varphi$ its Gauss map.

Let $T^3=\R^3/\Gamma$ be a 3-torus equipped with the standard flat metric. Let $x=(x_1,x_2,x_3)$ be the coordinate function on $\R^3$, which we also use as local flat coordinates on $T^3$. Now assume that $\Sigma\subset T^3$ is a closed (immersed) oriented surface of genus $g$ with induced metric and let $(u,v)$ be local isothermal coordinates on $\Sigma$. It is well-known that the (conformal class of) induced metric on $\Sigma$ determines a complex structure on $\Sigma$ and $z=u+iv$ is a local holomorphic coordinate. The metric on $\Sigma$ can be expressed as
\[\ud s^2=\rho(u,v)(\ud u^2+\ud v^2)=\rho(z,\bz)\ud z\ud\bz.\]
In other words,
\[\langle{\p x\over \p u},{\p x\over \p u}\rangle=\langle{\p x\over\p v},{\p x\over\p v}\rangle=\rho\textrm{\quad and \quad}\langle{\p x\over\p u},{\p x\over\p v}\rangle=0,\]
where $\langle\cdot,\cdot\rangle$ is the ambient metric on $T^3$. Moreover we can write down the K\"ahler form $\hat{\omega}$ as
\[\hat{\omega}=\rho(u,v)\ud u\wedge\ud v=\frac{i\rho(z,\bz)}{2}\ud z\wedge\ud\bz.\]
Locally, define
\[\phi=(\phi_1,\phi_2,\phi_3)={\p\over\p z}(x_1,x_2,x_3)=\frac{1}{2}({\p\over\p u}-i{\p\over\p v})(x_1,x_2,x_3).\]
The isothermal condition implies that
\[\phi_1^2+\phi_2^2+\phi_3^2=0.\]
It is a well-known result that $\Sigma$ is a minimal surface if and only if $\phi$ is holomorphic, or equivalently, $x=(x_1,x_2,x_3)$ is harmonic with respect to the Laplace-Beltrami operator
\[\Delta=\frac{1}{\rho}({\p^2\over\p u^2}+{\p^2\over\p v^2})=\frac{4}{\rho}{\p^2\over\p z\p\bz}\]
on $\Sigma$. Notice that for any smooth function $u$ on $\Sigma$, we have
\[2i\ddb u=\Delta u\cdot\hat{\omega}.\]

Now let us assume that $\Sigma$ is minimal and denote by $Q$ the Fermat quadric
\[Q=\{[\phi_1:\phi_2:\phi_3]\in\C\pp^2:\phi_1^2+\phi_2^2+\phi_3^2=0\}.\]
Therefore we get a map $\nu:\Sigma\to Q$ given by
\[z=u+iv\mapsto[\phi_1(z):\phi_2(z):\phi_3(z)],\]
which is a globally defined holomorphic map and it does not depend on the choice of local isothermal coordinates. This is known as the \emph{tangential Gauss map}.

A simple adjunction formula calculation indicates that $Q$ is biholomorphic to $\C\pp^1$. Let $\clo(1)$ be the positive generator of the Picard group of $Q$ and let $H$ be the hyperplane line bundle on $\C\pp^2$. It is easy to see that
\[H|_Q\cong\clo(2).\]
Moreover each $\phi_j$ serves as a section of $H$, which corresponds to a globally defined holomorphic 1-form $\mu_j=\phi_j\ud z$ on $\Sigma$. From this, we see that
\[\nu^*H\cong K_\Sigma,\]
in other words, $(\Sigma,\nu)$ is a vanishing spinorial pair.

In fact, an explicit identification $g:\C\pp^1\to Q$ is given by
\[[z_1:z_2]\mapsto [\phi_1:\phi_2:\phi_3]=[2z_1z_2:z_2^2-z_1^2:-i(z_1^2+z_2^2)].\]
Write $\zeta=z_2/z_1$, then
\[\zeta=\frac{\phi_2+i\phi_3}{\phi_1}=-\frac{\phi_1}{\phi_2-i\phi_3}.\]
Through the stereographic projection (\ref{stereographic}), it follows that the composition $\varphi=g^{-1}\circ\nu:\Sigma\to \C\pp^1=S^2$ is exactly the classical Gauss map defined by the unit normal vector field. Moreover, the pullback condition implies that $\varphi:\Sigma\to \C\pp^1$ is of degree $g-1$.

In fact, one can reconstruct the induced metric $\hat{\omega}$ from the holomorphic data $\mu_j=\phi_j\ud z$ through the so-called Weierstrass representation by the explicit formula
\begin{equation}\label{canonical}
\hat{\omega}=i(\mu_1\wedge\bar{\mu}_1+\mu_2\wedge\bar{\mu}_2+\mu_3\wedge\bar{\mu}_3).
\end{equation}

The classical Calabi-Gray case teaches us how to construct a canonical metric on generalized Calabi-Gray manifolds as follows.

Starting from a vanishing spinorial pair $(\Sigma,\varphi)$ and a fixed Fubini-Study metric $\omega_{\textrm{FS}}$ on $\ccc\pp^1$, one first chooses an orthonormal basis $\{\phi_j\}_{j=1}^3$ of holomorphic sections of $\clo(2)$ with respect to the $L^2$-inner product. Because of the pullback condition, the pullback sections $\mu_j:=\varphi^*\phi_j$ can be identified as holomorphic 1-forms on $\Sigma$, therefore we may define a canonical metric $\hat{\omega}$ on $\Sigma$ using the same recipe (\ref{canonical}). Notice that $\hat{\omega}$ does not depend on the choice of orthonormal basis and is uniquely defined up to scaling.

An alternative way to describe $\hat{\omega}$ is to identify $\ccc\pp^1$ with $Q$, we obtain a holomorphic map
\[g\circ\varphi:\Sigma\to Q\]
such that
\[(g\circ\varphi)^*(H|_Q)\cong K_\Sigma.\]
By pulling back $\phi_1$, $\phi_2$ and $\phi_3$ as sections of $H|_Q$, we get holomorphic 1-forms $\mu_1$, $\mu_2$ and $\mu_3$ on $\Sigma$ such that
\[\mu_1^2+\mu_2^2+\mu_3^2=0.\]
Again formula (\ref{canonical}) gives the canonical metric $\hat{\omega}$.

According to the work of Meeks \cite{meeks1990} and Traizet \cite{traizet2008}, for every $g\geq3$, there exist immersed and embedded minimal surfaces of genus $g$ in $T^3$. Conversely, it is an easy exercise in algebraic geometry that the existence of vanishing spinorial pair implies that the genus is at least three. Moreover, the moduli space of genus $g(\geq 3)$ curves admitting vanishing spinorial pairs form a divisor in the moduli space of curves \cite{harris1982, teixidor1987}. Therefore the generalized Calabi-Gray construction indeed gives more examples compared to the classical one.

A pair $(\Sigma,\varphi)$ satisfying the pullback condition is not so far from being minimal. It turns out that $\{\mu_j\}_{j=1}^3$ provides the Weierstrass data for the universal cover of $\Sigma$ to be minimally immersed into $\rr^3$ such that its Gauss map is given by $\varphi$ itself. In this sense, we may think of $(\Sigma,\varphi)$ with pullback condition satisfied as ``minimal surfaces in $T^3$'' without solving the period problem.

The simplest yet most useful example of vanishing spinorial pairs appears in the case $g=3$. As $S$ is of degree $g-1=2$ and it is basepoint-free, one concludes that $\Sigma$ must be a hyperelliptic curve and $\varphi:\Sigma\to\ccc\pp^1$ must be the hyperelliptic double covering, branching over 8 points on $\ccc\pp^1$. In fact, the moduli space of genus 3 curves with vanishing spinorial pair is exactly the hyperelliptic locus, which is parameterized by 8 distinct points on $\ccc\pp^1$. Moreover, for genus $3$ curves, $\{\mu_j\}_{j=1}^3$ forms a basis of holomorphic 1-forms on $\Sigma$, therefore by picking a reference point $p\in\Sigma$, integrating $\mu_j$ over cycles gives an embedding of $\Sigma$ into its Jacobian variety $J(\Sigma)$. It is not surprising that $\hat{\omega}$ is the induced metric on $\Sigma$ from the flat metric on $J(\Sigma)$. As $\Sigma$ is a complex submanifold of $J(\Sigma)$, by the curvature decreasing property, we know that the curvature of $\hat{\omega}$ is less or equal than the ambient curvature, therefore we conclude that the Gauss curvature $\kappa$ of $\Sigma$ is non-positive.

The fact that $\kappa\leq0$ holds for $\hat{\omega}$ constructed from any vanishing spinorial pair $(\Sigma,\varphi)$. One quick way to see this is that $(\Sigma,\hat{\omega})$ is locally a minimal surface in $\rr^3$, therefore its Gauss curvature is non-positive. A more detailed calculation reveals that
\[-\kappa\widehat{\omega}=i\pt\bpt\log\rho=\varphi^*\omega_{\textrm{FS}}\]
and
\[\|\ud\varphi\|^2_{\hat{\omega}}=-2\kappa,\]
therefore
$\kappa$ vanishes exactly at ramification points of $\varphi:\Sigma\to\ccc\pp^1$, and the vanishing order is twice of the ramification index minus two. In particular, for the hyperelliptic genus 3 curve case discussed above, $\kappa$ vanishes exactly at 8 hyperelliptic points with vanishing order 2 for each.

With the canonical metric $\hat{\omega}$ on $\Sigma$ understood, we can write down the product metric $\omega_0$ on $X=\Sigma\times M$ as
\[\omega_0=\hat{\omega}+\alpha\omega_I+\beta\omega_J+\gamma\omega_K,\]
which satisfies the balanced condition
\[\ud(\omega_0^2)=0.\]

\subsection*{Solving the Hull-Strominger System}

Another great thing about generalized Calabi-Gray manifolds is that the expression of the holomorphic volume $\Omega$ is explicit. As implicitly described in \cite{hitchin1987}, we may write
\[\Omega=\mu_1\wedge\omega_I+\mu_2\wedge\omega_J+\mu_3\wedge\omega_K.\]
With respect to the natural product metric $\omega_0$, we have
\[\|\Omega\|_{\omega_0}\equiv\textrm{const},\]
therefore the conformally balanced equation (\ref{cb})
\[\ud\left(\|\Omega\|_{\omega_0}\omega_0^2\right)=0\]
is satisfied.

Inspired by Fu-Yau \cite{fu2008}, we may deform the metric $\omega_0$ by adding a conformal factor $e^f$ depending only on $\Sigma$:
\begin{equation}\label{ansatz}
\omega_f=e^{2f}\widehat{\omega}+e^f(\alpha\omega_I+\beta\omega_J+\gamma\omega_K).
\end{equation}
For convenience of notation, we write
\[\omega'=\alpha\omega_I+\beta\omega_J+\gamma\omega_K.\]
An elementary calculation indicates that $\omega_f=e^{2f}\widehat{\omega}+e^f\omega'$ solves the conformally balanced equation (\ref{cb}) for arbitrary $f$.

The idea is very simple: we adopt the ansatz metric $\omega_f$ in (\ref{ansatz}) and try to find a suitable $f$ and a gauge bundle to solve the Hermitian-Yang-Mills equation (\ref{hym}) and the anomaly cancellation equation (\ref{ac}).

The calculation of the curvature term is rather complicated (see \cite{fei2017b, fei2017}). Here we only state the result:
\[\tr(R_f\wedge R_f)=i\pt\bpt\left(\frac{\|\ud\varphi\|_{\hat{\omega}}^2}{e^f}\omega'\right)+\tr(R'\wedge R'),\]
where $R'$ is the curvature form of the relative cotangent bundle $\Omega_{X/\Sigma}$ with respect to the metric induced from $\omega_0$. Therefore the $\tr(R'\wedge R')$-term can be cancelled by the curvature term from the gauge bundle, and the anomaly cancelation equation (\ref{ac}) reduces to
\begin{equation}\label{form1}
i\pt\bpt\left(\left(e^f+\frac{\alpha'\kappa}{2e^f}\right)\omega'\right)=0.
\end{equation}
This choice of gauge bundle makes sense since the relative cotangent bundle solves (\ref{hym}) automatically for arbitrary $\omega_f$. Roughly speaking this is because the fiber metrics of $p:X\to\Sigma$ are hyperk\"ahler, and we refer to \cite{fei2017b} for more details. For the case $M=T^4$, one finds that $\tr(R'\wedge R')=0$ though $R'\neq0$, so one may replace our choice of gauge bundle by any flat bundle. For the case $M=\textrm{K}3$, our choice of gauge bundle $E=\Omega_{X/\Sigma}$ has structure group $\U(2)$ instead of $\SU(2)$, however it is not hard to enhance it to an $\SU(N)$-bundle solving the Hull-Strominger system.

Notice that $u:=e^f+\dfrac{\alpha'\kappa}{2e^f}$ is a function depending only on $\Sigma$, and (\ref{form1}) is equivalent to
\begin{equation}
-\Delta u+2\kappa u=0.
\end{equation}

Therefore we have demonstrated that to solve the full Hull-Strominger system on generalized Calabi-Gray manifolds, we may use the ansatze (\ref{ansatz}) and the whole system reduces to a quadratic algebraic equation coupled with a linear elliptic PDE on $\Sigma$:
\begin{equation}\label{rac}\begin{cases}
&e^f+\dfrac{\alpha'\kappa}{2e^f}=u,\\
&\Delta u-2\kappa u=0.
\end{cases}\end{equation}

In particular, we obtain smooth solutions if and only if we can find a function $u$ in the kernel of the operator $-\Delta+2\kappa$ such that $u$ is positive at all ramification points of $\varphi$. Consequently the solvability of the Hull-Strominger system on $X$ is intimately related to spectral properties of the operator $-\Delta+2\kappa=-\Delta-\|\ud\varphi\|_{\hat{\omega}}^2$ on $\Sigma$. It is also clear that there are no solutions to the Hull-Strominger system under our ansatze if $\alpha'\leq 0$.

The operator $-\Delta+2\kappa=-\Delta-\|\ud\varphi\|_{\hat{\omega}}^2=:L_\varphi$ falls into the larger class of Sch\"odinger operators associated with holomorphic maps from Riemann surfaces to complex manifolds discussed in \cite{montiel1991}.

In the classical Calabi-Gray case, $\varphi$ is the Gauss map of a minimal surface $\Sigma$ in flat 3-space, hence $L_\varphi=-\Delta+2\kappa$ is the Jacobi operator, or the stability operator of the minimal surface, which originates from the second variation of the area functional.

Now let $\varphi:\Sigma\to\ccc\pp^1$ be a nonconstant holomorphic map. We would like to understand the space $\ker L_\varphi$. By the stereographic projection (\ref{stereographic}), we may think of $\varphi=(\alpha,\beta,\gamma)$ as a map into $\rr^3$. It is easy to check that $\alpha$, $\beta$ and $\gamma$ live in the kernel of $L_\varphi$, therefore
\[\dim\ker L_\varphi\geq 3.\]
It appears that in general we do not know how to compute $\ker L_\varphi$ and there are examples where $\dim L_\varphi$ is strictly greater than 3, though it is widely believed that $\dim\ker L_\varphi=3$ holds generically. For this reason, to solve the reduced the Hull-Strominger system (\ref{rac}) on generalized Calabi-Gray manifolds, we consider only the case $u\in\spann\{\alpha,\beta,\gamma\}$. In this scenario, the condition that $u$ is positive at all ramification points of $\varphi$ is equivalent to that all branched points of $\varphi$ on $\ccc\pp^1$ lie in an open hemisphere, which we abbreviate as the \emph{hemisphere condition}. Clearly if $(\Sigma,\varphi)$ satisfies the hemisphere condition and $\tau:\widetilde{\Sigma}\to\Sigma$ is an unramified covering of Riemann surfaces, then $(\widetilde{\Sigma},\varphi\circ\tau)$ also satisfies the hemisphere condition.

Given any vanishing spinorial pair $(\Sigma,\varphi)$, we may make the hemisphere condition holds by composing $\varphi$ with a suitable (hence infinitely many) automorphism of $\ccc\pp^1$, since there exist M\"obius transformations on $\ccc\pp^1$ pushing all points on the south hemisphere to the north pole. As a consequence, there exist vanishing spinorial pairs satisfying the hemisphere condition for every genus $g\geq3$. However, it seems unknown whether there are minimal surfaces in $T^3$ with its Gauss map satisfying the hemisphere condition. The hemisphere condition fails for the 5-dimensional family of triply periodic minimal surfaces constructed by Meeks \cite{meeks1990}.

In summary, we have sketched the proof of the following result:

\begin{thm}\cite{fei2017}\label{main}~\\
Let $(\Sigma,\varphi)$ be a vanishing spinorial pair with the hemisphere condition satisfied. Then we may construct explicit solutions to the Hull-Strominger system on the associated generalized Calabi-Gray manifold $X$. In fact, we get a family of such solutions of real dimension $\dim\ker L_\varphi\geq3$. As a consequence, for every genus $g\geq3$, there exist smooth solutions to the Hull-Strominger system on genus $g$ generalized Calabi-Gray manifolds. They have infinitely many distinct topological types and sets of Hodge numbers.
\end{thm}

It is worth mentioning that the hemisphere condition is also necessary when $\Sigma$ is of genus 3. In this case $\varphi:\Sigma\to\ccc\pp^1$ is the hyperelliptic covering, and we may use the hyperelliptic involution to help us.

\subsection*{A Discussion of Related Problems}

For any solution to the Hull-Strominger system, the uniqueness is always a hard question to answer. Ideally all the solutions are parameterized by some finite dimensional object such as cohomology classes. For example, a modified version of the Hull-Strominger system is proposed by Garcia-Fernadez et al. \cite{garciafernandez2018} where all the solutions are conjecturally parameterized by suitable Aeppli cohomology classes. However, as most, if not all, solutions known so far build on special ansatze, it is more reasonable to ask instead: are the solutions unique among the chosen ansatze?

For the Fu-Yau equation, the uniqueness is only established recently \cite{phong2018} in a suitable set, which is needed for ellipticity.

As for the generalized Calabi-Gray case, miracally the reduced system is automatically elliptic due to the sign of $\alpha'$ and $\kappa$ (this can be seen more clearly from the parabolic point of view in Section 4), we may expect to prove a better uniqueness result without further assumptions.

From the construction described in last subsection, our solution depends on a function $u\in\ker L_\varphi$ such that $u$ is positive at all ramification points of $\varphi:\Sigma\to\ccc\pp^1$. Suppose $\dim\ker L_\varphi=n\geq 3$ and we choose a basis $\{u_1,\dots,u_n\}$ of $\ker L_\varphi$. Let
\[D=\{u\in\spann\{u_1, u_2, \dots, u_n\}:u>0\textrm{ at ramification points of }\varphi\}.\]
The set $D$ is an open convex polyhedral cone in $\ker L_\varphi$. Consider the map $T:D\to\rr^n$ given by
\[T(u)=\left(\int_\Sigma e^fu_1\cdot\hat{\omega},\dots,\int_\Sigma e^fu_n\cdot\hat{\omega}\right),\]
where $e^f$ is determined from $u$ by the relation
\begin{equation} \label{e^f-to-u}
e^f=\frac{1}{2}\left(u+\sqrt{u^2-2\alpha'\kappa}\right)>0.
\end{equation}

It turns out that there exists a strictly convex function $F:D\to\rr$ such that $T(u)=\nabla F(u)$, where $F$ has the explicit form
\[F(u)=\frac{1}{2}\int_\Sigma \left(e^{2f}-\alpha'\kappa f\right)\hat{\omega}.\]
As a consequence, $T$ is a diffeomorphism from $D$ onto its image.

Since our ansatze solves the conformally balanced equation (\ref{cb}), the closed 4-form $\|\Omega\|_\omega\cdot\omega^2$ defines a de Rham cohomology class. In fact, this class is given by
\[\left[\|\Omega\|_{\omega_f}\cdot\omega_f^2\right]=2\sqrt{2}\left([\vol_M]+ [\omega_I][e^f\alpha\cdot\hat{\omega}]+[\omega_J][e^f\beta\cdot\hat{\omega}]+ [\omega_K][e^f\gamma\cdot\hat{\omega}]\right)\in H^4(X;\rr)\]
under the K\"unneth isomorphism
\[H^*(X;\rr)\cong H^*(M;\rr)\otimes H^*(\Sigma;\rr).\]
Therefore the conformally balanced class is parameterized by three numbers
\[\int_\Sigma e^f\alpha\cdot\widehat{\omega},\quad \int_\Sigma e^f\beta\cdot\widehat{\omega}, \quad\textrm{ and }\quad \int_\Sigma e^f\gamma\cdot\widehat{\omega}.\]

If $\dim\ker L_\varphi=3$, which is expected to hold generically, then $\ker L_\varphi$ is spanned by $\alpha$, $\beta$ and $\gamma$, hence $T$ can be regard as a ``period map'' which maps a solution to the Hull-Strominger system to its associated de Rham conformally balanced class. Since $T$ is a diffeomorphism onto its image, we obtain
\begin{cor}\cite{fei2017}~\\
Assuming the ansatze (\ref{ansatz}) and $\dim\ker L_\varphi=3$, then there is at most one solution to the Hull-Strominger system in each de Rham cohomology class in $H^4(X;\rr)$.
\end{cor}
In this case, the period map $T$ has a potential function $F:D\to\rr$. A natural question is that does $F$ has a critical point in $D$? In other words, does the period map $T$ hit the origin? Our conjecture is that $F$ has no critical points in $D$, which is partially supported by the explicit example to be given in next subsection.

If $M=T^4$, it makes sense to talk about a different kind of uniqueness. In this case the corresponding generalized Calabi-Gray manifolds have a large symmetry group: all translations on $T^4$ are holomorphic. Furthermore in this case, we may take the gauge bundle to be flat, so we only need to solve for the Hermitian metric $\omega$ on $X$ from the system:
\begin{equation}\label{treduce}
\begin{cases}&\ud\left(\|\Omega\|_\omega\cdot\omega^2\right)=0,\\
&i\pt\bpt\omega=\dfrac{\alpha'}{4}\tr(R\wedge R).
\end{cases}
\end{equation}
The question is: are the solutions we found the only $T^4$-invariant solutions to the reduced system (\ref{treduce})?

Another interesting problem is about geometric consequences of the hemisphere condition. If we consider a general holomorphic map $\varphi:\Sigma\to\ccc\pp^1$ and its associated elliptic operator $L_\varphi$, we have
\begin{thm}\cite{tysk1987, grigoryan2004}
\[C\cdot\deg\varphi\leq\ind~L_\varphi\leq 7.68183\cdot\deg\varphi.\]
\end{thm}
Here $C$ is a universal positive constant and $\ind~L_\varphi$ is the number of negative eigenvalues of the operator $L_\varphi$. However, if the hemisphere condition is imposed, then we have the better estimate \cite{fei2017}:
\[\ind~L_\varphi\geq\deg\varphi,\]
which does not hold if the hemisphere condition fails. It is very interesting to explore other geometric implications of the hemisphere condition and its possible relationship with stability conditions.

After solving the Hull-Strominger system mathematically on generalized Calabi-Gray manifolds, one naturally asks whether one can build a feasible physical theory. A first attempt was recently made by Chen-Pantev-Sharpe \cite{chen2018}. We look forward to further developments in this direction.

\subsection*{An Explicit Example}

For expository purpose, let us end this section with an explicit example of generalized Calabi-Gray manifold of genus 3.

Let $a$ be a positive real number. Let $\Sigma$ be the genus 3 hyperelliptic Riemann surface defined by the equation
\[w^2=\zeta^8-a^8.\]
The set $\{(w,\zeta)\in\C^2:w^2=\zeta^8-a^8\}$ is an affine chart of $\Sigma$, which consists of $\Sigma$ with two points removed. By implicit function theorem, $\zeta$ is a local holomorphic coordinate when $w\neq0$ and $w$ is a local holomorphic coordinate when $\zeta\neq0$.

The other affine chart can be identified as $\{(v,\eta)\in\C^2:v^2=1-a^8\eta^8\}$,
where the transformation is given by
\[\eta=\frac{1}{\zeta},\quad v=\frac{w}{\zeta^4}.\]
Moreover, the hyperelliptic degree 2 map is
\[\varphi:(w,\zeta)\mapsto \zeta\in\pp^1,\]
and the hyperelliptic involution is $\iota:(w,\zeta)\mapsto(-w,\zeta)$.

Now we can write down $\phi_1$, $\phi_2$ and $\phi_3$ as sections of $\clo(2)=T\pp^1$ explicitly as
\[\begin{split}\phi_1&=\zeta\frac{\p}{\p\zeta}=-\eta\frac{\p}{\p\eta},\\ \phi_2&=\frac{1}{2}(\zeta^2-1)\frac{\p}{\p\zeta}=\frac{1}{2}(\eta^2-1)\frac{\p}{\p\eta},\\ \phi_3&=-\frac{i}{2}(\zeta^2+1)\frac{\p}{\p\zeta}=\frac{i}{2}(\eta^2+1)\frac{\p}{\p\eta}.\end{split}\]

As $\varphi^*\left(\dfrac{\p}{\p\zeta}\right)$ is a globally defined holomorphic 1-form of $\Sigma$ with 2 double zeroes at fibers over $\infty$, or $(v,\eta)=(\pm1,0)$, we can write down
\[\varphi^*\left(\dfrac{\p}{\p\zeta}\right)=\frac{\eta^2}{v}\ud\eta=-\frac{\ud v}{4a^8\eta^5}=-\frac{\ud\zeta}{w}=-\frac{\ud w}{4\zeta^7}.\]
Therefore we can write
\[\begin{split}\mu_1&=-\frac{\zeta}{w}\ud\zeta,\\ \mu_2&=-\frac{\zeta^2-1}{2w}\ud\zeta,\\ \mu_3&=\frac{i(\zeta^2+1)}{2w}\ud\zeta.\end{split}\]
Therefore we can write down the metric
\[\widehat{\omega}=\frac{i(|\zeta|^2+1)^2}{2|w|^2}\ud\zeta\wedge\ud\bzeta=\frac{i(|\zeta|^2+1)^2}{2|\zeta^8-a^8|} \ud\zeta\wedge\ud\bzeta = \frac{i(|\zeta|^2+1)^2}{32|w|^{14}}\ud w\wedge\ud\bw.\]

The eight branched points are
\[\zeta=ae^{k\pi i/4}\]
for $k=0,1,\dots,7$, which corresponds to
\[(\alpha,\beta,\gamma)=\left(\frac{1-a^2}{1+a^2},\frac{2a\cos\frac{k\pi}{4}}{1+a^2},\frac{2a\sin \frac{k\pi}{4}}{1+a^2}\right).\]
For simplicity let us assume that $0<a<1$ so all the branched points satisfy $\alpha>0$, i.e. they all lie on an open hemisphere.

In this case one can compute that
\[-\kappa=\frac{4|\zeta^8-a^8|}{(|\zeta|^2+1)^4},\]
hence
\[\max(-\kappa)=4\]
with equality holds at $\zeta=\infty$.

For $u=t_1\alpha+t_2\beta+t_3\gamma$, the hemisphere condition is
\[(1-a^2)t_1+2at_2\cos\frac{k\pi}{4}+2at_3\sin\frac{k\pi}{4}>0\] for $k=0,1,2,\dots,7$. In particular we may take $t_1>0$ and $t_2=t_3=0$.

Now consider $u=t\alpha=t\cdot\dfrac{1-|\zeta|^2}{1+|\zeta|^2}$, then
\[e^f(t)=\frac{1}{2}\left(u+\sqrt{u^2-2\alpha'\kappa}\right)=\frac{1}{2(1+|\zeta|^2)}\left(t(1-|\zeta|^2) +\sqrt{t^2(1-|\zeta|^2)^2+8\alpha'\frac{|\zeta^8-a^8|}{(1+|\zeta|^2)^2}}\right).\]
By symmetry it is not hard to show that
\[\int e^f\beta=\int e^f\gamma=0.\]
Moreover, we have
\[\begin{split}&\frac{1}{2}\int e^f\alpha=\int_{\pp^1}\frac{i(1-|\zeta|^2)}{2|\zeta^8-a^8|}\left(t(1-|\zeta|^2) +\sqrt{t^2(1-|\zeta|^2)^2+8\alpha'\frac{|\zeta^8-a^8|}{(1+|\zeta|^2)^2}}\right)\ud\zeta\wedge\ud\bzeta\\ =&\int_0^{2\pi}\int_0^\infty\frac{1-r^2}{|r^8e^{8i\theta}-a^8|}\left(t(1-r^2) +\sqrt{t^2(1-r^2)^2+8\alpha'\frac{|r^8e^{8i\theta}-a^8|}{(1+r^2)^2}}\right)r\ud r\ud\theta\\ =&\int_0^\pi\int_0^\infty\frac{1-s}{|s^4e^{i\theta}-a^8|}\left(t(1-s) +\sqrt{t^2(1-s)^2+8\alpha'\frac{|s^4e^{i\theta}-a^8|}{(1+s)^2}}\right)\ud s\ud\theta\\ =&\int_0^\pi\int_0^1t(1-s)^2\left(\frac{1}{|s^4e^{i\theta}-a^8|}+\frac{1}{|e^{i\theta}-s^4a^8|}\right)\ud s\ud\theta\\ +&\int_0^\pi\int_0^1\frac{1-s}{1+s}\left(\frac{\sqrt{t^2(1-s^2)^2+8\alpha'|s^4e^{i\theta}-a^8|}}
{|s^4e^{i\theta}-a^8|}-\frac{\sqrt{t^2(1-s^2)^2+8\alpha'|e^{i\theta}-s^4a^8|}}{|e^{i\theta}-s^4a^8|}\right)\ud s\ud\theta.\end{split}\]

Notice that
\[|e^{i\theta}-s^4a^8|^2-|s^4e^{i\theta}-a^8|^2=(1-s^8)(1-a^{16})>0,\]
we conclude that $\int e^f\alpha$ is always positive.

It follows that in this particular example we are considering, the potential function $F$ does not have any critical point in $D$ hence our conjecture in last subsection holds.

Based on all these information, we may easily proceed to write down the explicit solution metric to the Hull-Strominger system on the associated generalized Calabi-Gray manifold as well.

\section{Anomaly Flow and Its Reduction on Riemann Surfaces}

As we have seen in previous sections, all the known solutions to the Hull-Strominger rely heavily on a clever choice of ansatze and then a reduction of the original system to a simpler one. However if one wishes to consider the general situation, we can no longer choose an ansatze. Alternatively, we can only let the geometry and analysis choose the path to solution for us. Based on this idea, Phong-Picard-Zhang \cite{phong2018b} initiated the Anomaly flow program as a parabolic approach to the Hull-Strominger system which we shall describe.

Mathematically, the idea is to evolve the conformally balanced metric by the Green-Schwarz term
\[\pt_t\left(\|\Omega\|_\omega\cdot\omega^2\right)=i\pt\bpt\omega-\frac{\alpha'}{4}(\tr(R\wedge R)-\tr(F\wedge F)),\]
while simultaneously flow the Hermitian metric on the gauge bundle by the Donaldson heat flow
\[H^{-1}\pt_tH = -\Lambda F.\]
It is obvious that the Anomaly flow preserves the conformally balanced equation (\ref{cb}) and its stationary points are exactly the solutions to the Hull-Strominger system.

Given knowledge of the bundle metric, Phong-Picard-Zhang \cite{phong2018e} showed that this flow is equivalent to evolving the Hermitian metric by the second Chern-Ricci curvature rescaled by the dilaton function plus corrections from curvature square term:
\[\pt_t\omega(t)=-\frac{\rho^{(2)}(t)}{2\|\Omega\|_{\omega(t)}}+\textrm{ correction terms},\]
where $\rho^{(2)}(t)$ is the second Chern-Ricci form associated to the evolving metric $\omega(t)$.

In this regard, we may think of the Anomaly flow as a generalization of the K\"ahler-Ricci flow into the non-K\"ahler setting. In their paper, Phong-Picard-Zhang were able to establish the short-time existence of the flow when the initial $\|\alpha'\textrm{Rm}\|$ is small, as well as the long-time existence (for $\alpha'=0$) under the assumption of non-collapsing and that the torsion, the first derivative of torsion and the curvature tensors do not blow up \cite{phong2018e}, i.e.,
\[\|Rm\|^2+\|T\|^4+\|DT\|^2<\infty.\]

From the point view of the Hull-Strominger system, the Anomaly flow is very natural in the sense that every known reduction in the elliptic side has its parabolic version. For instance, Phong-Picard-Zhang \cite{phong2018d} gave a uniform treatment of the Fu-Yau solution in regardless of sign of $\alpha'$ and investigated the behavior of the Anomaly flow on complex Lie groups \cite{phong2017b}. In a very recent preprint \cite{phong2018c}, Phong-Picard-Zhang showed that on K\"ahler Calabi-Yau manifolds, the Anomaly flow with $\alpha'=0$ preserves the conformally K\"ahler structure and converges to the Ricci-flat K\"ahler metric.

The same phenomenon also occurs in the setting of generalized Calabi-Gray manifolds. As observed in \cite{fei2017d}, the Anomaly flow under the ansatze (\ref{ansatz}) reduces to an evolution equation on the Riemann surface $\Sigma$:
\begin{equation}\label{pcalabi}
\pt_t(e^f)=\Delta\left(e^f+\frac{1}{2}\alpha'\kappa e^{-f}\right)-2\kappa\left(e^f+\frac{1}{2}\alpha'\kappa e^{-f}\right).
\end{equation}
In terms of the function $u:=e^f+\frac{1}{2}\alpha'\kappa e^{-f}$, the above equation is equivalent to
\[\pt_t u=\left(1-\frac{1}{2}\alpha'\kappa e^{-2f}\right)(\Delta u-2\kappa u),\]
which is automatically parabolic due to the fact that $\alpha'>0$ and $\kappa\leq0$. As a consequence we get the short-time existence for free.

Regarding (\ref{pcalabi}), we have the following results:
\begin{thm}\cite{fei2017d}
\begin{enumerate}
\item The flow can be continued as long as $e^f$ is bounded from below by a positive number,
\item If $\|e^f\|_{L^1(\Sigma,\hat{\omega})}$ is sufficiently small, then the flow develops finite-time singularity and we have an estimate of the maximal existence time. This can be interpreted as that if the initial area of the Riemann surface is small, then its area keeps decreasing. In general, it does not shrink uniformly like the Ricci flow, singularities appear at certain points when the Riemann surface still has positive area.
\item If $e^f\geq C$ initially for some constant $C$, then we have long-time existence and the flow is uniformly parabolic. In this case $e^f$ grows exponentially and the growth rate is given by the first eigenvalue of the Jacobi operator $\Delta-2\kappa$, for which one can estimate by the dimension reduction method \cite{fei2018}. After suitable normalization, $e^f$ converges to the first eigenfunction of $\Delta-2\kappa$. In this case, the condition $\|\alpha'Rm\|\ll1$ is preserved and in fact it decays exponentially. Similarly the quantity $\|Rm\|^2+\|T\|^4+\|DT\|^2$ also decays exponentially. Geometrically this can be interpreted as that generalized Calabi-Gray manifolds collapse their hyperk\"ahlr fibers under the Anomaly flow and converge to the Riemann surface base in the Gromov-Hausdorff sense.
\end{enumerate}
\end{thm}

The above-mentioned results leave the region of medium initial data, where we have the stationary points of the flow. This case is very interesting since one would expect to detect the ``hemisphere condition'', which is an obstruction to the existence of solutions. Naively one should expect singularities corresponding to the failure of hemisphere condition. For these singularities, the flow can be continued after reparametrization (surgery).

In summary, the Hull-Strominger system and the Anomaly flow lie at the crossroad of string theory, algebraic geometry and geometric analysis. Its reduction on Riemann surfaces connects the theory of algebraic curves, minimal surfaces, parabolic flows and potentially to complex function theory and degenerate equations. What we know currently is only a tip of iceberg, there are a great number of exciting problems awaiting us to investigate!
\\

\noindent\textbf{Acknowledgements: }The author would like to thank his coauthors Zhijie Huang, Sebastien Picard and Shing-Tung Yau for the wonderful collaborative experience. The author is also indebted to Duong H. Phong, Baosen Wu, Chenglong Yu and Xiangwen Zhang for many inspiring discussions leading to his work and this paper. Finally the author thanks Xiangwen Zhang for carefully reading the manuscript and giving many suggestions on writing.

\bibliographystyle{plain}

\bibliography{C:/Users/Piojo/Dropbox/Documents/Source}

\bigskip
Department of Mathematics, Columbia University, New York, NY 10027, USA

\smallskip
tfei@math.columbia.edu

\end{document}